\documentclass[10pt]{amsart}
\usepackage{amsmath, latexsym, amssymb, color}
\usepackage[all]{xy}
\def\a{{\alpha}}
\def\b{{\beta}}

\def\e{{\epsilon}}
\def\w{{\omega}}
\def\g{{\gamma}}

\newcommand\res{\operatorname{res}}
\def\dummy{{\mbox{-}}}

\def\F{{\mathcal{F}}}

\def\BC{{\mathbb{C}}}
\def\BZ{{\mathbb{Z}}}
\def\BQ{{\mathbb{Q}}}
\def\Hom{{\operatorname{Hom}}}
\def\End{{\operatorname{End}}}

\def\bidu{^{\vee\vee}}
\newcommand\Unit{U}
\newcommand\Counit{C}
\newcommand\triv{{\operatorname{triv}}}

\newcommand\C[1]{{#1\mbox{-\bf{mod}}_{\operatorname{\mathsf{fin}}}}}
\newcommand\CM[1]{{#1\mbox{-\bf{mod}}}}
\newtheorem{thm}{Theorem}[section]

\newtheorem{defn}[thm]{Definition}
\newtheorem{prop}[thm]{Proposition}
\newtheorem{remark}[thm]{Remark}

\newtheorem{lem}[thm]{Lemma}

\newtheorem{example}[thm]{Example}

\newcommand\ptr{\operatorname{\underline{ptr}}}
\newcommand\ptrl{\ptr^\ell}
\newcommand\ptrr{\ptr^r}
\newcommand\FS{\operatorname{\underline{FS}}}
\newcommand\id{\operatorname{id}}
\renewcommand\o{{\otimes}}
\newcommand\ro{\otimes}
\newcommand\Tr{\operatorname{Tr}}

\newcommand\piv{j}
\newcommand\dual[1]{#1^\vee}

\newcommand\inv{^{-1}}
\DeclareMathOperator\ev{\operatorname{ev}}
\DeclareMathOperator\db{\operatorname{db}}
\newcommand\CC{\mathcal C}

\newcommand\FF{\mathcal F}
\newcommand\bidual[1]{#1^{\vee\vee}}
\newcommand\du{^{\vee}}

\newcommand{\ou}[1]{\mathrel{\mathop{\otimes}_{#1}}}

\newcommand\so[1]{^{(#1)}}

\newcommand\hh[1]{\widehat{#1}}
\title[Frobenius-Schur Indicators for Quasi-Hopf  Algebras]
{Central Invariants and Higher Indicators for Semisimple Quasi-Hopf  Algebras}
\author{Siu-Hung Ng}
\address{Department of Mathematics, Iowa State University, Ames, IA 50011, USA}
\email{rng@iastate.edu}
\thanks{The first author is supported by the NSA grant number H98230-05-1-0020.}
\author{Peter Schauenburg}
\address{Mathematisches Institut der Universit\"at M\"unchen,
Theresienstr.\ 39, 80333 M\"unchen, Germany}
\email{schauenburg@math.lmu.de}
\thanks{The second author is supported by a DFG Heisenberg fellowship.}
\begin{document}
\begin{abstract}
  In this paper, we define the higher Frobenius-Schur (FS-)indicators for finite-dimensional
  modules $V$ of a semisimple quasi-Hopf algebra $H$ via the categorical counterpart developed in \cite{NS05}.
  We prove that this definition of higher FS-indicators coincides with the higher indicators
  introduced by Kashina, Sommerh\"auser, and Zhu when $H$ is a Hopf algebra. We also obtain a sequence of canonical central
  elements of $H$, which is invariant under gauge transformations, whose values, when evaluated by
  the character of an $H$-module $V$, are the higher Frobenius-Schur indicators of $V$.
  As an application, we show that FS-indicators are sufficient to distinguish the four
  gauge equivalence classes of semisimple quasi-Hopf algebras of dimension eight
  corresponding to the four fusion
  categories with certain fusion rules classified by Tambara and Yamagami.
  Three of these categories correspond to well-known Hopf algebras, and we explicitly construct
  a quasi-Hopf algebra corresponding to the fourth.
   We also derive explicit formulae for FS-indicators for some quasi-Hopf algebras
   associated to group cocycles.
\end{abstract}
\maketitle
\section*{Introduction}
The notion of (degree 2) Frobenius-Schur indicator $\nu_2(V)$ of an
irreducible representation $V$ of a finite group $G$ has been
generalized to simple modules of semisimple Hopf algebras by
Linchenko and Montgomery \cite{LM00}, to certain $C^*$-fusion
categories by Fuchs, Ganchev, Szlach\'anyi, and Vescerny\'es
\cite{FGSV99}, and to simple modules of semisimple quasi-Hopf
algebras by Mason and the first author \cite{MN04}. A more general
version of the Frobenius-Schur Theorem holds for the simple modules
of semisimple Hopf algebras or even quasi-Hopf algebras. In
particular, the Frobenius-Schur indicator of a simple module is
non-zero if, and only if,
the simple module is self-dual, and its value can only be  0, 1 or -1. \\

In proving that $\pm 1$ are the only possible non-zero values for
the Frobenius-Schur indicator of a simple module of a semisimple
quasi-Hopf algebra $H$ over $\BC$ \cite{MN04}, the fact that $\C{H}$
is pivotal, proved by Etingof, Nikshych, and Ostrik \cite{ENO}, has
been used. Based on the pivotal structure, the second author
\cite{Sch04} later introduced a categorical definition of degree 2
Frobenius-Schur indicators and
gave a different proof of the Frobenius-Schur Theorem for quasi-Hopf algebras.\\

The higher indicators of irreducible representations of a finite
group do not have a direct interpretation as the degree 2 indicators
(cf. \cite{IsaacsBook}). The $n$-th Frobenius-Schur indicator of a
finite-dimensional module $V$ with character $\chi$ of a semisimple
Hopf algebra $H$ has been formally defined by Kashina,
Sommerh\"auser, and Zhu \cite{KSZ} as the value of $\chi$ at the
$n$-th Sweedler power of the normalized integral of $H$. It has been
shown in their paper that these indicators carry rich information on
$H$, as well as its module category. Moreover, the values of the
$n$-th indicators are cyclotomic integers in
$\BQ_n$.\\

In the paper \cite{NS05}, the authors introduced the notion of
higher Frobenius-Schur indicators of an object $V$ in a $\BC$-linear
pivotal category $\CC$. These indicators are invariants of the
pivotal category. If $V$ is simple and $\nu_n(V) \ne 0$, then the
dual object $V\du$ of $V$ occurs in $V^{\o (n-1)}$. Again, the
$n$-th indicator of any object in $\CC$ is a cyclotomic integer in
$\BQ_n$. In addition, if $\CC$ is a {\em pseudo-unitary fusion
category over $\BC$}, then the higher indicators will be preserved
by any $\BC$-linear monoidal equivalence from $\CC$ to any other
pseudo-unitary fusion
category over $\BC$.\\

In this paper, we will adopt the categorical definition of higher
Frobenius-Schur indicators given in \cite{NS05} to define the higher
Frobenius-Schur indicators of a finite-dimensional module $V$ of a
semisimple quasi-Hopf algebra $H$ over $\BC$ with respect to the
canonical pivotal (or pseudo-unitary) structure of the category
$\C{H}$. We show that for any positive integer $n$, there exists a
canonical central element $\mu_n(H)$ such that $n$-th
Frobenius-Schur indicator $\nu_n(V)$ of $V$ is given by
$$
\nu_n(V)=\chi(\mu_n(H))
$$
where $\chi$ is the character of the $H$-module $V$. Moreover,
$\mu_n(H)$ is invariant under gauge transformations, and it is
independent of the choice of antipode of $H$. This formula implies
that our definition of higher indicators coincides with
the one introduced in \cite{KSZ} when $H$ is a Hopf algebra.\\

The organization of the paper is as follows: we cover some basic
definitions and facts about quasi-Hopf algebras $H$, including some
important elements of $H \o H$ and identities introduced by Drinfeld
\cite{Drin90},  Hausser, and Nill \cite{HNQA}, in Section \ref{s2}.
In Section \ref{s:H-mod}, we prove that two finite-dimensional Hopf
algebras over a field $k$ are gauge equivalent if, and only if,
their module categories are $k$-linear monoidally equivalent. In
addition, if $H$ is semisimple and $k=\BC$, then $\C{H}$ is a
spherical fusion category with respect to a canonical pivotal
structure. Moreover, the canonical pivotal structure of $\C{H}$ is
preserved by any $\BC$-linear monoidal equivalence from $\C{H}$ to
$\C{K}$ for some quasi-Hopf algebra $K$ over $\BC$. In Sections
\ref{s:FSI} and \ref{s:FSE}, we define the $(n,r)$-th
Frobenius-Schur indicators $\nu_{n,r}(V)$ of a finite-dimensional
$H$-module $V$. We determine the central elements $\mu_n(H)$ whose
action on the $H$-modules gives the $n$-th Frobenius-Schur
endomorphisms introduced in \cite{NS05}. The element $\mu_n(H)$ is
an invariant under gauge transformations on $H$ and
$\nu_{n}(V)=\chi(\mu_n(H))$ where $\chi$ is the character of $V$. As
an example, we derive a formula for the higher indicators for the
semisimple quasi-Hopf algebra obtained from a semisimple Hopf
algebra with a central group-like element of order 2 in Section
\ref{s:HCG}. In Section 6, we use this formula to show that
Frobenius-Schur indicators suffice to identify and distinguish the
four gauge equivalence classes of non-commutative semisimple
quasi-Hopf algebras of dimension 8 whose fusion rules (or $\mathcal
K(H)$) contain an abelian group isomorphic to $\mathbb
Z_2\times\mathbb Z_2$. The corresponding categories were classified
by Tambara and Yamagami. Finally, in Section \ref{s:DG}, we obtain
formulae for the higher FS-indicators for the dual of the group
algebra of a finite group $G$ with quasi-Hopf algebra structure
determined by a group 3-cocycle, and for its double, known as
the twisted double $D^{\w}(G)$.\\

\subsection*{Acknowledgements}
The authors are indebted to Susan Montgomery for her interest in
their paper and extended correspondence on the results in Sections
\ref{s:HCG} and \ref{ea}. The first author would like to thank
Penn State University and IH\'{E}S for their hospitality during
the preparation of this manuscript. The second author thanks the
DFG for support by a Heisenberg fellowship. He also thanks the
Institute of Mathematics of the University of Tsukuba and Akira
Masuoka for their hospitality.

\section{Preliminaries and Notations}\label{s2}
In this section, we recall the definition of quasi-Hopf algebras,
some properties described in \cite{Drin90} and \cite{Kassel}, and
some interesting results obtained in \cite{HNQA}, \cite{HN991},
\cite{HN992}. In the
sequel, we will use the notation introduced in this section. \\

A {\em quasi-bialgebra} over a field $k$ is a quadruple $(H, \Delta,
\e, \phi)$, in which $H$ is an algebra over $k$, $\Delta\colon H
\rightarrow H \otimes H $ and $\e \colon H \rightarrow k$ are
algebra maps, and $\phi \in H^{\o 3}$ is the {\em associator}. Here
our convention of associator $\phi$ is given by the equation
$$
\phi (\Delta \otimes \id)\Delta(h)= (\id\otimes \Delta
)\Delta(h)\phi\,.
$$
A quasi-bialgebra $(H, \Delta, \e,\phi)$ is called a {\em quasi-Hopf
algebra} if there exist an anti-algebra automorphism $S$ of $H$ and
elements $\a, \b \in H$ such that for all elements $h \in H$, we
have
\begin{equation}\label{2.8}
  S(h_{(1)})\a h_{(2)} =\e(h)\a, \quad  h_{(1)}\b S(h_{(2)})
  =\e(h)\b\,,
\end{equation}
\begin{equation}\label{2.9}
  \phi^{(1)} \b S(\phi^{(2)})\a \phi^{(3)} =1, \quad
   S(\phi^{(-1)}) \a \phi^{(-2)}\b S(\phi^{(-3)}) =1\,,
\end{equation}
where $\phi=\phi^{(1)} \o \phi^{(2)}\o \phi^{(3)}$, $\phi^{-1}=
\phi^{(-1)} \otimes \phi^{(-2)}\otimes \phi^{(-3)}$ and
$\Delta(h)=h_{(1)}\otimes  h_{(2)}$. In the above equations, the
summation notations of the tensors have been suppressed. For
simplicity, we will continue to do so in the sequel. We will simply
write $H$ for a quasi-bialgebra $(H, \Delta, \e,\phi)$ or a
quasi-Hopf algebra $(H, \Delta, \e,\phi, \a, \b , S)$.\\

The module category $\CM{H}$ of the quasi-bialgebra $H$ is a
$k$-linear monoidal category, or simply a {\em tensor category}. If
$H$ is a quasi-Hopf algebra, then the tensor category $\C{H}$ of all
finite-dimensional $H$-modules is {\em rigid}, i.e. $\C{H}$ admits
both left and right duality. Given any $V \in \C{H}$ with the left
$H$-module structure given by $\rho\colon H \rightarrow \End_k(V)$,
the left dual $(V\du, \ev, \db)$ of $V$ is defined as follows:
\begin{enumerate}
\item  $V\du = \Hom_k(V,k)$ with the $H$-action given by $h \mapsto
\rho(S(h))^*$,
\item $\ev\colon V\du \o V \rightarrow k$ and $\db\colon k \rightarrow V
\o V\du$  defined by
$$
\ev(f \o v)= f(\a v)\quad \mbox{and} \quad \db(1)= \sum_i \b v_i \o
v^i\,,
$$
where $\{v_i\}$ is a basis for $V$ and $\{v^i\}$ the
corresponding dual basis.
\end{enumerate}
 Similarly, one
can define ${\du V}=\Hom_k(V,k)$ with the left $H$-action given by
$h \mapsto \rho(S\inv(h))^*$ and the linear maps $\ev'\colon V \o
{\du V} \rightarrow k$ and  $\db'\colon k \rightarrow {\du V} \o V$
by
$$
\ev'(v \o f)= f(S\inv(\a) v), \quad \mbox{and} \quad \db'(1)= \sum_i
v^i \o S\inv(\b) v_i\,.
$$
Then $({\du V}, \ev', \db')$ defines a right dual of $V$ (cf.
\cite{Drin90} and \cite{Kassel} for the
details).\\

Following \cite{Kassel}, a {\em gauge transformation} on a
quasi-bialgebra $H=(H, \Delta, \e, \phi)$ is an invertible element
$F$ of $H \otimes H$ such that
$$
(\e \otimes \id)(F)=1 =(\id \otimes \e)(F)\,.
$$
Using a gauge transformation $F$ on $H$, one can define an algebra
map $\Delta^F\colon H \rightarrow H \otimes H$ by
\begin{equation}\label{2.6}
  \Delta^F(h) = F \Delta(h) F^{-1}
\end{equation}
for any $h \in H$, and an invertible element $\phi^F$ of $H \otimes
H \otimes H$ by
\begin{equation}\label{2.7}
  \phi^F=(1 \otimes F) (\id \otimes \Delta)(F) \phi (\Delta \otimes
  \id)(F^{-1})(F^{-1} \otimes 1)\,.
\end{equation}
Then $H^F=(H, \Delta^F, \e, \phi^F)$ is  also a quasi-bialgebra. In
addition, if $H=(H, \Delta, \e, \phi, \a, \b, S)$ is a quasi-Hopf
algebra, then so is $H^F=(H, \Delta^F, \e, \phi^F, \a^F, \b^F, S)$,
where
$$ \a^F = \sum_i S(d_i) \a e_i
\quad\mbox{and}\quad \b^F = \sum_i f_i \b S(g_i)
$$
with $F=\sum_i f_i \otimes g_i$ and $F^{-1} =\sum_i d_i \otimes
e_i$.\\

Two quasi-bialgebras $A$ and $B$ are said to be {\em gauge
equivalent} if there exists a gauge transformation $F$ on $B$ such
that $A$ and $B^F$ are isomorphic as quasi-bialgebras. Let $\sigma
\colon A \rightarrow B^F$ be such a quasi-bialgebra isomorphism.
Then the functor ${_\sigma (\dummy)}\colon \CM{B} \rightarrow
\CM{A}$, with ${_\sigma V}$ the left $A$-module with the underlying
space $V$ and the left $A$-action given by
\begin{equation}\label{eq:sigmafunctor}
 a\cdot v = \sigma(a)v \quad (a \in A, v \in V),
\end{equation}
and ${_\sigma f}=f$ for any map $f$ in $\CM{B}$, is a $k$-linear
equivalence. Let $\xi \colon {_\sigma V} \otimes  {_\sigma
W}\rightarrow {_\sigma (V \otimes W)}  $ be the linear map
$$
V \otimes W \xrightarrow{F \cdot} V \otimes  W
$$
for any $V , W \in \CM{B}$. Then $({_\sigma (\dummy)}, \xi, \id)$ is
a $k$-linear monoidal equivalence, or simply a {\em tensor
equivalence}, from
$\CM{B}$ to $\CM{A}$ (cf. \cite{Kassel}). \\

In \cite{HNQA}, \cite{HN991} and \cite{HN992}, Frank Hausser and
Florian Nill introduced some interesting elements in $H \otimes H$
for any arbitrary quasi-Hopf algebra $H=(H, \Delta, \e,\phi, \a, \b
, S)$ in the course of studying the corresponding theories of
quantum double, integral and the fundamental theorem for quasi-Hopf
algebras. These elements of $H \otimes H$ are given by
\begin{eqnarray}
\label{R} & q_R  = \phi^{(1)} \otimes S^{-1} (\a
\phi^{(3)})\phi^{(2)} \,, &\quad p_R= \phi^{(-1)} \otimes
\phi^{(-2)} \b
S(\phi^{(-3)})\,, \\
\label{L} & q_L = S(\phi^{(-1)}) \a \phi^{(-2)} \otimes
\phi^{(-3)}\,, &\quad  p_L =  \phi^{(2)}S^{-1}(\phi^{(1)}\b) \otimes
\phi^{(3)} \,.
\end{eqnarray}
The elements $q_L$ and $p_L$ also occurred in \cite{Drin90}. One can
show (cf. \cite{HNQA}) that they obey the relations (for all $a \in
H$)
\begin{equation} \label{eq:qRpR}
 (a\otimes 1)\, q_R =   (1 \otimes S^{-1}(a_{(2)}))\,q_R\, \Delta(a_{(1)}),\quad
  p_R\, (a \otimes 1) =   \Delta(a_{(1)})\, p_R \, (1 \otimes  S(a_{(2)})),
  \end{equation}
\begin{equation} \label{eq:qLpL}
(1 \otimes a)\, q_L =  (S(a_{(1)})\otimes 1)\, q_L \, \Delta(a_{(2)}), \quad
p_L\, (1 \otimes a) =    \Delta(a_{(2)})\,p_L\,
(S^{-1}(a_{(1)})\otimes 1)\,,
\end{equation}
where $\Delta(a)= a_{(1)} \otimes a_{(2)}$. Suppressing the
summation symbol and indices again, we write $q_R = q^{(1)}_R
\otimes q^{(2)}_R$, etc. These elements also satisfy the identities
(cf. \cite{HNQA}):
\begin{equation} \label{eq:qRpR2}
\Delta(q_R^{(1)}) \, p_R\, (1\otimes S(q_R^{(2)})) = (1\otimes S^{-1}(p_R^{(2)}))\, q_R \, \Delta(p_R^{(1)})
= 1\otimes 1,
\end{equation}
\begin{equation} \label{eq:qLpL2}
\Delta (q_L^{(2)}) \, p_L \, (S^{-1}(q_L^{(1)}) \otimes 1)=(S(p_L^{(1)}) \otimes 1)\, q_L \, \Delta(p_L^{(2)})
 = 1\otimes 1.
\end{equation}

\section{Module Categories of Quasi-Hopf Algebras}\label{s:H-mod}
In this section, we recall the canonical pivotal structure of the
module categories of finite-dimensional semisimple quasi-Hopf
algebras over $\BC$ and some properties of these tensor categories.
We also prove that two finite-dimensional quasi-Hopf algebras are
gauge equivalent if,
and only if, their module categories are tensor equivalent.\\

It is well-known that if $H$ and $K$ are gauge equivalent
quasi-bialgebras, then $\CM{H}$ and $\CM{K}$ are equivalent tensor
categories (cf. \cite{Kassel}). The converse for Hopf algebras was
proved in \cite{S96}. The quasi-bialgebra case was proved in
\cite[section 6]{EG02}. Here we give a more straightforward proof
for the case of finite-dimensional quasi-Hopf algebras. \\

\begin{lem}\label{l:tensor}
  Let $H$ be a quasi-Hopf algebra over a field $k$, $R$ a $k$-algebra  and $V$ an
  $H$-$R$-bimodule.
Then $\theta\colon H \otimes V \rightarrow H \otimes {_\circ V}$ given by
   $$
   \theta(h \otimes v) = q_R^{(1)} h_{(1)} \otimes S(q_R^{(2)} h_{(2)})v,
  $$

  for any $h \in H$ and $v \in V$, is a natural $H$-$R$-bimodule isomorphism, where
  ${_\circ V}$ denotes the trivial $H$-module with the underlying space $V$ and
  the right $R$-actions on $H \otimes V$ and $H \otimes {_\circ V}$ are induced
  by the right $R$-action on $V$.
\end{lem}
\begin{proof}
  It follows directly from \eqref{eq:qRpR} that  $\theta$ is natural bimodule homomorphism.
  Consider the $H$-module map
  $\bar{\theta}\colon H \otimes {_\circ V} \rightarrow
   H \otimes  V$ given by
   $$
   \bar{\theta}(h \otimes v) = h_{(1)} p_R^{(1)}  \otimes h_{(2)} p_R^{(2)} v\,.
   $$
   Using \eqref{eq:qRpR} and \eqref{eq:qRpR2},
   one can easily verify that
   $\theta\bar{\theta}=\bar{\theta}\theta=\id$. Hence, the result follows.
\end{proof}

\begin{thm}\label{t:gauge_equivalence}
  Let $H$ be a finite-dimensional quasi-Hopf algebra, and $B$ a
  quasi-bialgebra over a field $k$. If
  $\CM{H}$ and $\CM{B}$ are tensor equivalent (in particular if $B$ is finite-dimensional,
  and $\C H$ and $\C B$ are tensor equivalent), then $B$ is gauge equivalent to $H$
  as quasi-bialgebras.
\end{thm}
\begin{proof}
  Let $(\F, \xi, \xi_0)$ be a tensor equivalence from $\CM{H}$ to $\CM{B}$.
  By the Morita Theorems, we may assume that  $\FF= T\ou H \!-$
  for some $T\in B\mbox{-}\CM{H}$ such that $T_H$ is a
  progenerator, and in particular, $T$ is finite-dimensional.
  Moreover, the algebra homomorphism $\sigma'\colon B \rightarrow \End_H(T_H)$
  given by the $B$-module structure of $T$, is an isomorphism. To show that
  $H \cong B$ as algebras, it suffices to prove that $T\cong H$ as right
  $H$-modules. By Lemma \ref{l:tensor}, there exists an $H$-$H$-bimodule isomorphism
  $\theta\colon H \o H \rightarrow H \o {_\circ H}$. Thus, we have
  \begin{multline*}
    T\o T \cong
       \FF(H)\o\FF(H)
       \stackrel{\xi}{\cong}\FF(H\o H)
       \stackrel{\FF(\theta)}{\cong}\FF( H\o {_\circ H})
       = T \o_H H \o {_\circ H}
       \cong T\o {_\circ H}
  \end{multline*}
  as left $B$-modules. Obviously, all the above unlabeled isomorphism are $B$-$H$-bimodule isomorphisms.
  Since $\theta$ is an $H$-$H$-bimodule isomorphism, $\FF(\theta)$ is a $B$-$H$-bimodule
  isomorphism.
  The naturality of
  $\xi$ implies that $\xi \colon \FF(H) \o \FF(H)\rightarrow \FF(H\o H)$ is also
  a $B$-$H\o H$-bimodule isomorphism. Thus, we have
  $$
  T \o T \cong T \o {_\circ H}
  $$
  as  $B$-$H$-bimodules and hence $T \cong H$ as right $H$-modules by the Krull-Schmidt Theorem.
  Let $\sigma\colon B \rightarrow H$ be the composition map
  $$
  B \xrightarrow{\sigma'} \End_H(T_H) \cong  \End_H(H_H) \cong H\,.
  $$
  As in \eqref{eq:sigmafunctor}, the algebra map $\sigma$ induces a
  $k$-linear equivalence $_\sigma(\dummy): \C{H} \rightarrow \C{B}$
  and $T \cong {_\sigma H}$ as $B$-$H$-bimodules. The following
  $B$-module isomorphisms
  $$
  T \o_H V \cong {_\sigma H} \o_H V \cong {_\sigma V}
  $$
  are natural in $V$ and hence $_\sigma (\dummy)$ is $k$-linearly equivalent to $\FF$.
  Therefore, one may assume
  $\FF(\dummy)={_\sigma(\dummy)}$. Note that $(\FF, z \xi, z\inv \xi_0)$ is also a
  tensor equivalence for any non-zero scalar $z$. One may further
  assume $\xi_0=\id_k$ and so we have $\e_H \sigma = \e_B$. Let
  $$
  F' =  \xi_{H,H} (1 \o 1)\quad \mbox{and}\quad F= \xi_{H,H}\inv (1 \o 1)\,.
  $$
  Since $\xi_{H,H}$ is a $B$-$H\o H$-bimodule
  isomorphism,  we have
  $$
  \xi_{H,H}(u \o v) = \xi_{H,H}(1 \o 1)(u \o v)=F'(u \o v)
  $$
  for all $u$, $v \in H$, and
   $$
  F'(\sigma \otimes \sigma)(\Delta_B (b)) = \xi(b \cdot (1 \o 1))
  =b \cdot  \xi(1 \o 1) = \Delta_H(\sigma(b))F'\,.
  $$
  for all  $b \in B$. By naturality again, $\xi_{X,Y} = F'\cdot \,$ for any
  free $H$-modules $X$ and $Y$. Moreover, we have
  $$
  FF'=\xi\inv(1 \o 1)(F')=\xi\inv(F')=\xi\inv(\xi(1 \o 1))=1 \o
  1
  $$
  and, similarly,  $F'F=1 \o 1$. Therefore, $F$ is invertible in $H \o H$.
  Finally, by the commutativity of
  the diagrams
  $$
  \xymatrix{
  & {_\sigma H}&\\
  {_\sigma H} \o {_\sigma H}  \ar[rd]_-{
  \e \o \id}\ar[rr]^{\xi} \ar[ru]^-{\id \o \e}
  && {_\sigma (H \o H)}\ar[ld]^-{\e \o \id}\ar[lu]_-{\id \o \e}\,\,,\\
  &{_\sigma H}&
  }
  \quad
  \xymatrix{{_\sigma ((H \otimes  H) \otimes  H)}\ar[r]^-{_\sigma (\phi_H\,\cdot)}
     &
    {_\sigma (H \otimes  (H\otimes  H))}
    \\
    {_\sigma (H \otimes  H) \otimes  _\sigma H}\ar[u]_-{\xi}
     & {_\sigma H} \otimes  {_\sigma (H\otimes  H)}
    \ar[u]_-{\xi } \\
    ({_\sigma H} \otimes  {_\sigma H}) \otimes {_\sigma H}
  \ar[r]^-{\phi_B \cdot} \ar[u]_-{\xi \otimes \id} &
  {_\sigma H} \otimes ( {_\sigma H} \otimes {_\sigma H})
  \ar[u]_-{\id \otimes \xi }\,\,,
  }
  $$
  we obtain
  $$
  1 =(\id \otimes \e_H  )(F) = (\e_H \otimes \id)(F)\,,
  $$
  and
  $$
  (\sigma \otimes \sigma \otimes \sigma)(\phi_B)=
  (1 \otimes F)(\id \o \Delta) (F)\phi_H (\Delta \otimes \id)(F^{-1}) (F^{-1} \otimes
  1)\,.
  $$
  Therefore, $F$ is a gauge transformation on $H$ and $\sigma\colon B \rightarrow H^F$ is quasi-bialgebra
  isomorphism.
\end{proof}

Let us further assume that $H$ is a finite-dimensional semisimple
quasi-Hopf algebra over $\BC$. It follows from \cite[Section 8]{ENO}
that there exists a unique pivotal structure on $\C{H}$, i.e. an
isomorphism of tensor functors $j: Id \rightarrow (\dummy)\bidu$,
such that the {\em left pivotal dimension} $\ptrl(\id_V)$ of any
finite-dimensional $H$-module $V$ is identical to its usual
dimension $\dim(V)$, where
$$
\ptrl(f):=\left(\BC\xrightarrow{\db} V \o V\du \xrightarrow{j_V \o
\id_{V\du}} V\bidu \o V\du \xrightarrow{\ev} \BC\right)
$$
for any $f \in \End_H(V)$.
\begin{remark}\label{r:spherical}
Let $V \in \C{H}$ and $ V = \bigoplus_i X_i $ a decomposition  of
$V$ as a direct sum of simple $H$-modules $X_i$. Let $\iota_i \colon
X_i \rightarrow V$ and $\pi_i \colon V \rightarrow X_i$ be the
embeddings and projections associated with the decomposition. For
any $f \in \End_H(V)$, we have $\pi_i f \iota_i = f_i \id_{X_i}$ for
some scalar $f_i \in \BC$ by Schur's Lemma. Therefore,
$$
\ptrl(f) = \sum_i \ptrl(f \iota_i \pi_i) = \sum_i \ptrl(\pi_i f \iota_i) =
\sum_i f_i \ptrl(\id_{X_i})= \sum_i f_i \dim(X_i)
$$
which is identical to the usual trace of $f$. In particular,
$$
\ptrl(f) = \ptrl(f\du)= \ptrr(f)\,.
$$
Therefore, $\C{H}$ is a (non-strict) spherical fusion category over
$\BC$. This observation has been proved, in a more general context,
by M\"{u}ger \cite[Lemma 2.8]{MugerI}.
\end{remark}

The pivotal structure of $\C{H}$ could be explicitly described in
terms of the {\em trace element} $g_H$ of $H$ (cf. \cite{MN04}). In
the sequel, we will consider $\C{H}$ as a spherical fusion category
with respect to the pivotal structure $j$ described above. Such
spherical fusion category is also called pseudo-unitary in
\cite{ENO}. This pivotal structure will be automatically preserved
by any tensor equivalence to any other pseudo-unitary fusion
category over $\BC$.

\begin{prop}\label{p:pivotal_preservation}
 Let $H$, $K$ be gauge equivalent finite-dimensional semisimple quasi-Hopf
 algebras
over $\BC$. Then every tensor equivalence from $\C{H}$ to $\C{K}$
preserves their underlying pivotal structures.
\end{prop}
\begin{proof}
  Since $\C{H}$ and $\C{K}$ are pseudo-unitary, the statement follows
  immediately from \cite[Corollary 6.2]{NS05}.
\end{proof}
\section{Frobenius-Schur indicators for semisimple quasi-Hopf
algebras}\label{s:FSI} In this section, we recall the definition of
Frobenius-Schur indicators of an object in a linear pivotal category
(cf. \cite{NS05} for the details). We then give a definition of
higher Frobenius-Schur indicators for any finite-dimensional
representations of a semisimple quasi-Hopf algebra $H$ over $\BC$
using the canonical pivotal structure of $\C{H}$ described in
Section \ref{s:H-mod}. It follows from \cite{NS05} that these
indicators are invariants of the tensor category
$\C{H}$ and hence gauge invariants of $H$.\\

Let $\CC$ be a finite $k$-linear pivotal category; that is a
$k$-linear rigid monoidal category with a pivotal structure $j:Id
\rightarrow (-)\bidu$ such that $\CC(V,W)$ is a finite-dimensional
$k$-linear space for all $V, W \in \CC$. We denote by $V^{\ro n}$
the $n$-fold tensor power of an object $V\in\CC$ with rightmost
parentheses; thus $V^{\ro 0}=I$, the neutral object of $\CC$, and
$V^{\ro(n+1)}=V\o V^{\ro n}$. There is a unique isomorphism
$$\Phi^{(n)}\colon V^{\ro(n-1)}\o V\to V^{\ro n}$$ composed of
instances of the associativity isomorphisms $\Phi$; explicitly
$\Phi^{(1)}$ is the identity, and
$$\Phi^{(n+1)}=\left(
  (V\o V^{\ro(n-1)})\o V
  \xrightarrow{\Phi}V\o(V^{\ro(n-1)}\o V)
  \xrightarrow{V\o\Phi^{(n)}}
  V^{\ro(n+1)}\right).$$

For $V,W\in\CC$, define $A\colon \CC(I, V \o W) \rightarrow
\CC(V\du, W)$ and  $T_{VW}\colon\CC(\dual V,W)\to\CC(\dual W,V)$ by
$$
A(f)=\left(V\du \xrightarrow{\dual V\o f}\dual V\o (V\o
  W)\xrightarrow{\Phi\inv}(\dual V\o V)\o W\xrightarrow{\ev\o
  W}W\right),\\
$$
$$T_{VW}(f)=(\dual W\xrightarrow{\dual f}\bidual
V\xrightarrow{\piv_V\inv}V),$$ and put
\begin{gather*}E_{VW}=\left(\CC(I,V\o W)\xrightarrow{A}
  \CC(\dual V,W)\xrightarrow{T_{VW}}
  \CC(\dual W,V)\xrightarrow{A\inv}
  \CC(I,W\o V)\right),
  \\E^{(n)}_V=\left(\CC(I,V^{\ro n})\xrightarrow{E_{V,V^{\ro(n-1)}}}
                  \CC(I,V^{\ro{(n-1)}}\o
                  V)\xrightarrow{\CC(I,\Phi^{(n)})}
                  \CC(I,V^{\ro n})\right)\,.
\end{gather*}
Note that $E_V^{(1)}=\id_{\CC(I,V)}$ as $I\bidu=I$. Following
\cite{NS05}, for any positive integers $n, r$, the $(n, r)$-th
\textbf{Frobenius-Schur indicator} of $V$ is the scalar
$$
\nu_{n, r}(V) = \Tr\left(\left(E_V^{(n)}\right)^r\right)\,.
$$
We will call $\nu_{n}(V):=\nu_{n,1}(V)$ the $n$-th
\textbf{Frobenius-Schur indicator} of $V$. Now, we can define the
Frobenius-Schur indicators for the representations of a semisimple
quasi-Hopf algebra over $\BC$.
\begin{defn}\label{d:FSI}
  Let $H$ be a semisimple quasi-Hopf algebra over $\BC$ and let $\CC$ be the spherical category
  $\C{H}$ with respect to the pivotal structure described at the end of Section \ref{s:H-mod}.
  For any $V \in \CC$,
  we call $\nu_{n,r}(V)$ the $(n,r)$-th Frobenius-Schur indicator of $V$ and call
  $\nu_n(V)$ the $n$-th Frobenius-Schur indicator of $V$.
\end{defn}
\begin{prop}\label{p:invarianceI}
  Let $H$, $K$ be gauge equivalent finite-dimensional semisimple quasi-Hopf algebras over
  $\BC$ via the gauge transformation $F$ on $H$ and the quasi-bialgebra isomorphism
  $\sigma\colon K \rightarrow H^F$.
  For any tensor equivalence $\F$ from $\C{H}$ to $\C{K}$,
  $$
  \nu_{n,r}(V)  = \nu_{n,r}(\F(V))
  $$
  for any $V \in \C{H}$ and positive integers $n, r$. In particular,
  $$
  \nu_{n,r}(V)  = \nu_{n,r}({_\sigma V})\,.
  $$
\end{prop}
\begin{proof}
  The result follows directly from Proposition \ref{p:pivotal_preservation} and \cite[Corollary 4.4]{NS05}.
\end{proof}
\begin{remark}\label{r:invarianceI}
  {\rm
  Proposition \ref{p:invarianceI} implies that the $(n,r)$-th Frobenius-Schur indicators are
  gauge invariants of $H$.
  }
\end{remark}

\begin{remark}
  Let $H$ be a semisimple Hopf algebra over $\BC$.
  By a well-known result of Larson and Radford \cite{LaRa88}, the antipode of $H$ is an involution, and so
  the canonical pivotal structure of $\C{H}$ is given by the
natural isomorphism $j: V \rightarrow V\bidu$ of $\BC$-linear
spaces. If one identifies $\Hom_H(\BC,V)$ with the invariant space
$V^H$ for any $V \in \C{H}$, then
$$
E_V^{(n)}\left(\sum u_1 \o \cdots \o u_n\right) = \sum u_2 \o \cdots
\o u_n \o u_1
$$
for any $\sum u_1 \o \cdots \o u_n \in \left(V^{\o n}\right)^H$. By
\cite[Corollary 2.3]{KSZ}, $\nu_n(V)=\Tr(E_V^{(n)})$ is identical to
the $n$-th indicator of $V$ defined by Kashina, Sommerh\"{a}use and
Zhu \cite{KSZ}. Therefore, definition of higher indicators given in
Definition \ref{d:FSI} is indeed a generalization of the higher
indicators for Hopf algebras.
\end{remark}

\section{Frobenius-Schur Endomorphisms--Central Gauge Invariants}   \label{s:FSE}
In \cite{NS05}, we have defined the $n$-th Frobenius-Schur
endomorphism of an object in a semisimple pivotal monoidal category,
and related it to the Frobenius-Schur indicators. In the category
$\C{H}$ for a semisimple quasi-Hopf algebra $H$, the Frobenius-Schur
endomorphism is given by multiplication with a central element
$\mu_n(H)\in H$. In this section we will obtain an explicit formula
for this element $\mu_n(H)$ in terms of the quasi-Hopf algebra
structure and the normalized integral of $H$. In the case of an
ordinary Hopf algebra, this formula simplifies to the $n$-th
Sweedler power of the integral, so that $\chi(\mu_n(H))$ specializes
to the original definition of $\nu_n(V)$ in \cite{KSZ}.\\

Let $k$ be a field and $\CC$ a finite $k$-linear semisimple pivotal
category with pivotal structure $j\colon Id \rightarrow (-)\bidu$
and neutral object $I$ such that
$$
\CC(X,X)\cong \CC(I,I) \cong k
$$
for all simple objects $X$ of $\CC$. By \cite{NS05}, there exists a
unique natural isomorphism $\tau_{VT}\colon V\o T\to T\o V$ for any
$I$-isotypical object $T$ and for any $V\in\CC$ such that
$\tau_{VI}=\id_V$. One can define the $n$-th \textbf{Frobenius-Schur
endomorphism} of $V$ as the
  composition
  \begin{equation}\label{eq:FSE}
  \begin{aligned}
  \FS^{(n)}_V=\biggl(V\xrightarrow{\Unit}& Y\o(V^{\o(n-1)}\o V)
     \xrightarrow{Y\o\Phi^{(n)}}Y\o V^{\o n} \xrightarrow{Y\o\pi}\\
        & Y\o(V^{\o n})^\triv\xrightarrow{\tau}(V^{\o n})^\triv\o Y
     \xrightarrow{\iota\o Y}V^{\o n}\o Y
     \xrightarrow{\Counit}V\biggr),
  \end{aligned}
  \end{equation}
  where
  \begin{gather*}
    \Unit=\left(V\xrightarrow{\db'\o V}(Y\o V^{\o(n-1)})\o V
      \xrightarrow{\Phi}Y\o(V^{\o(n-1)}\o V)\right),\\
    \Counit=\left((V\o V^{\o(n-1)})\o
    Y\xrightarrow{\Phi}V\o(V^{\o(n-1)}\o
    Y)\xrightarrow{V\o\ev'}V\right),
  \end{gather*}
  $((V^{\o n})^\triv,\iota,\pi)$ is the $I$-isotypical
    component of $V^{\o n}$, and $(Y, \ev', \db')$ is a right dual of $V^{\o n}$.
    Moreover, by \cite[Theorem 7.6]{NS05},
\begin{equation}\label{eq:FSI-ptr}
\nu_n(V)= \ptrl(\FS^{(n)}_V)\,,
\end{equation}
where
$$
\ptrl(f)\colon=\left(I\xrightarrow{\db}V\o V\du\xrightarrow{f\o
V\du}V\o V\du\xrightarrow{j_V\o V\du}V\bidu\o
V\du\xrightarrow{\ev}I\right)
$$
for any $f \in \CC(V,V)$.
 In the above equation, the identification $\CC(I,I)=k$ has been used.\\

Now, let $\CC=\C{H}$ for some finite-dimensional semisimple
quasi-Hopf algebra $H=(H, \Delta, \e, \phi, \a, \b, S)$ over $\BC$.
As described in Section \ref{s2}, for any $V \in \CC$,  $({\du V},
\ev', \db')$ defines a right dual of $V$. Thus the maps $U$ and $C$
in the definition of Frobenius-Schur endomorphism can be expressed
in terms of $q$ and $p$ as follows: Let $\{u_i\}$ be a basis for
$V^{\o (n-1)}$ and $\{u^i\}$ its dual basis for $\du\left(V^{\o
(n-1)}\right)$. For any $x \in V$, $u \in V^{\o n}$ and $f \in
{\du\left(V^{\o (n-1)}\right)}$, we have
$$
U(x)=u^i \o p_L^{(1)}u_i \o p_L^{(2)}x,\quad \mbox{and}\quad
C(x \o u \o f)=q_R^{(1)}f(q_R^{(2)}u)x\,.
$$
By \cite{HNQA}, $H$ admits a unique normalized two-sided integral
$\Lambda$; that is the two-sided integral of $H$ such that
$\e(\Lambda)=1$. Then the trivial isotypical component $V^\triv$ of
a finite-dimensional $H$-module $V$ is given by $\Lambda V$ and
$\pi(x)=\Lambda x$ defines a retraction of the inclusion map
$\iota\colon V^\triv \rightarrow V$.\\

Let us define
$$
\Delta^{(0)}=\e, \quad \Delta^{(1)}=\id_H, \quad
\Delta^{(2)}=\Delta, \quad \phi_1=1_H, \quad \phi_2=1_H \o 1_H,
$$
and recursively
\begin{equation}\label{eq:phi_n}
\Delta^{(n+1)}=(\id_H \o \Delta)\Delta^{(n)}, \quad \phi_{n+1}
=(1 \o \phi_n) (\phi^{(1)} \o \Delta^{(n-1)}(\phi^{(2)})\o \phi^{(3)})
\end{equation}
for any positive integer $n \ge 2$. Then $\Phi^{(n)}\colon V^{\o
(n-1)} \o V \rightarrow V^{\o n}$ is  given by the action of
$\phi_n$ on $V^{\o n}$. For any $a \in H$, we will suppress the
summation notation and simply write
$$
a_{(1)} \o \cdots \o a_{(n)}
$$
for $\Delta^{(n)}(a)$. In this notation, for any $V \in \CC$, the
action of an element $a$ of $H$ on $V^{\o n}$ is given by
$$
a\cdot (x_1 \o \cdots \o x_n) = a_{(1)}x_1 \o \cdots \o a_{(n)}x_n
$$
for  $x_1 \o \cdots \o x_n \in V^{\o n}$.\\

Now we can derive a formula for the $n$-th Frobenius-Schur
endomorphism $\FS_V^{(n)}$ in $\CC$. By \eqref{eq:FSE}, we obtain
\begin{equation*}
\begin{aligned}
\FS^{(n)}_V(x)& =  \sum_{i_1, \cdots, i_{n-1}}
q_R^{(1)}\Lambda_{(1)}\phi_n^{(1)}p_{L, (1)}^{(1)} x_{i_1}
\langle q_{R, (1)}^{(2)}\Lambda_{(2)}\phi_n^{(3)}p_{L, (2)}^{(1)} x_{i_2}, x^{i_1}\rangle \\
& \cdots
\langle q_{R, (n-2)}^{(2)}\Lambda_{(n-1)}\phi_n^{(n-1)}p_{L, (n-1)}^{(1)}x_{i_{n-1}},  x^{i_{n-2}}\rangle
\cdot \langle q_{R, (n-1)}^{(2)}\Lambda_{(n)}\phi_n^{(n)}p_L^{(2)}x , x^{i_{n-1}}\rangle
\end{aligned}
\end{equation*}
where $\{x_i\}$ is a basis for $V$ and $\{x^i\}$ is its dual basis
for $V^*$. Since $\sum_i x_i \langle v, x^i \rangle = v$ for all $v
\in V$, one can simplify $\FS^{(n)}_V(x)$ as follows:
\begin{equation} \label{eq:FSE-mu}
\begin{aligned}
\FS^{(n)}_V& (x)  = \left(q_R^{(1)}\Lambda_{(1)}\phi_n^{(1)}p_{L,
(1)}^{(1)}\right)\cdot
                  \left(q_{R, (1)}^{(2)}\Lambda_{(2)}\phi_n^{(3)}p_{L, (2)}^{(1)}\right) \\
&\quad \quad \quad  \cdots \left(q_{R,
(n-2)}^{(2)}\Lambda_{(n-1)}\phi_n^{(n-1)}p_{L,
(n-1)}^{(1)}\right)\cdot
\left(q_{R, (n-1)}^{(2)}\Lambda_{(n)}\phi_n^{(n)}p_L^{(2)}\right)x\\
& = m\left((q_R^{(1)} \o \Delta^{(n-1)}(q_R^{(2)}))\cdot
\Delta^{(n)}(\Lambda)\cdot \phi_n \cdot (\Delta^{(n-1)}(p_L^{(1)})
\o p_L^{(2)})\right)x\,,
\end{aligned}
\end{equation}
where $m$ is the multiplication on $H$. Let us define
\begin{equation}\label{eq:mu}
\mu_n(H):= m\left((q_R^{(1)} \o \Delta^{(n-1)}(q_R^{(2)}))\cdot
\Delta^{(n)}(\Lambda)\cdot \phi_n \cdot (\Delta^{(n-1)}(p_L^{(1)})
\o p_L^{(2)})\right)
\end{equation}
for any integer $n \ge 1$. The following theorem shows that
$\mu_n(H)$ is a central gauge invariant of $H$.
\begin{thm}\label{t:1}
  Let $H=(H, \Delta, \e, \phi, \a, \b, S)$ be a finite-dimensional semisimple quasi-Hopf algebra over $\BC$
  and $n$ a positive integer. The element $\mu_n(H)$ is in the center of $H$ and it is invariant under
  gauge transformations on $H$. Moreover, for any $V \in \C{H}$,
  $$
  \nu_n(V)=\chi(\mu_n(H))
  $$
  where $\chi$ is the character of $V$. In addition, if both $\a$ and $\b$ are invertible elements of $H$,
  then the element $\mu_n(H)$ is  given by
  $$
  \mu_n(H) =m( \Delta^{(n)}(\Lambda)\phi_n) (\b \a)\inv = (\b \a)\inv m( \Delta^{(n)}(\Lambda)\phi_n)
  $$
  where $\Lambda$ is the normalized two-sided integral of $H$.
\end{thm}
\begin{proof}
  Since $\FS_H^{(n)}$ is an $H$-module map, the equality \eqref{eq:FSE-mu} implies that
  $\mu_n(H)$ must lie in the center of $H$. It follows from Remark \ref{r:spherical} and
  \eqref{eq:FSI-ptr} that
  $$
  \nu_n(V)=\ptrl(\FS_V^{(n)}) = \Tr(\FS_V^{(n)})=\chi(\mu_n(H))
  $$
  where $\Tr(\FS_V^{(n)})$ denotes the usual trace of the linear operator $\FS_V^{(n)}$. \\

  Let $K$ be
  another semisimple quasi-Hopf algebra over $\BC$ such that $K$ is gauge equivalent to $H$ via the
  gauge transformation $F$ on $H$ and the quasi-bialgebra isomorphism $\sigma\colon K  \rightarrow H^F$.
  Then $({_\sigma(\dummy)},F\cdot\,, \id)$ is a tensor equivalence from $\C{H}$ to $\C{K}$.
  By \cite[Lemma 7.3]{NS05}, the functor ${_\sigma(\dummy)}$ preserves Frobenius-Schur endomorphisms. Therefore,
  $$
  \sigma(\mu_n(K))x=\mu_n(K)\cdot x = \FS_K^{(n)}(x) =\FS_{\,_\sigma H}^{(n)}(x)= {_\sigma(\FS_H^{(n)})}(x)
  = \mu_n(H)x
  $$
  for all $x \in {_\sigma H}$. Therefore, $\sigma(\mu_n(K))=\mu_n(H)$. In particular, $\mu_n(H^F)=\mu_n(H)$.\\

  By \cite[Lemma 3.1]{MN04}, we have the equations
  \begin{equation}\label{eq:q_Lambda_p}
  \b q_X^{(1)}\Lambda_{(1)} \o q_X^{(2)}\Lambda_{(2)} = \Delta(\Lambda)=
   \Lambda_{(1)}p_Y^{(1)} \o \Lambda_{(2)}p_Y^{(2)} \a
  \end{equation}
  for any $X,Y \in \{R,L\}$.
  Following from \eqref{eq:q_Lambda_p} and the equation
  \begin{equation}\label{eq:phi_nLambda}
  \phi_n (\Delta^{(n-1)}  \o
  \id)\Delta(\Lambda)=\Delta^{(n)}(\Lambda)\phi_n\,,
  \end{equation}
   we have
   $$
   (\b q_R^{(1)} \o \Delta^{(n-1)}(q_R^{(2)}))\cdot \Delta^{(n)}(\Lambda)\cdot \phi_n
\cdot (\Delta^{(n-1)}(p_L^{(1)}) \o p_L^{(2)}\a) =
\Delta^{(n)}(\Lambda)\phi_n\,.
   $$
  Since $\mu_n(H)$ is in the center of $H$,
  $$
  \b\a\mu_n(H)=\mu_n(H)\b\a = \b \mu_n(H)\a=m\left(\Delta^{(n)}(\Lambda)\phi_n\right)\,.
  $$
  In addition, if $\b$ and $\a$ are invertible elements of $H$, the last statement follows.
\end{proof}

  \begin{remark} Since two semisimple quasi-Hopf algebras $H$, $H'$ with
  identical quasi-bialgebra structures but different antipodes are
  gauge equivalent via the gauge transformation $1\o 1$ and the
  quasi-bialgebra isomorphism $\id_H$, Theorem \ref{t:1} implies that
  $\mu_n(H)=\mu_n(H')$. Therefore, $\mu_n(H)$ is independent of the
  choice of antipode of $H$.
\end{remark}

\begin{remark}
  Let $H$ be a semisimple Hopf algebra over $\BC$ and $V$ a finite-dimensional $H$-module with character $\chi$.
  By Theorem \ref{t:1}, $\nu_n(V)=\chi(m(\Delta^{(n)}(\Lambda)))=\chi(\Lambda^{[n]})
  $ which coincides with the $n$-th Frobenius-Schur indicator of $V$
  for Hopf algebras introduced in \cite{KSZ}.
\end{remark}
Equation \eqref{eq:mu} suggests that the other combinations of the
$q$'s and $p$'s may not give the same element. Indeed, $\mu_n(H)$
can be expressed in any of the four combinations of those $q$'s and
$p$'s. We need the following observation for the proof.

\begin{lem}\label{product_trick}
  Let $n$ be a positive integer and let $t$ be an element of $H^{\otimes
  (n+1)}$. For any $G \in H^{\otimes n}$,
  $$
   m\left((1 \otimes G)t  \right) = m\left(t(G\otimes 1)\right)
   $$
   where $m$ is the multiplication on $H$.
\end{lem}
\begin{proof}
The statement can be easily verified by direct computation.
\end{proof}
\begin{prop}
 Let $H=(H, \Delta, \e, \phi, \a, \b, S)$ be a finite-dimensional semisimple quasi-Hopf algebra over $\BC$
  and $n$ a positive integer. Then
  $$
  \mu_n(H)= m\left((q_X^{(1)} \o \Delta^{(n-1)}(q_X^{(2)}))\cdot \Delta^{(n)}(\Lambda)\cdot \phi_n
\cdot (\Delta^{(n-1)}(p_Y^{(1)}) \o p_Y^{(2)})\right)
  $$
  for any $X, Y \in \{R, L\}$.
\end{prop}
\begin{proof} Let
$$
T_{X,Y}=m\left((q_X^{(1)} \o \Delta^{(n-1)}(q_X^{(2)}))\cdot \Delta^{(n)}(\Lambda)\cdot \phi_n
\cdot (\Delta^{(n-1)}(p_Y^{(1)}) \o p_Y^{(2)})\right)
$$
for any $X,Y \in \{R, L\}$. By definition, $T_{R,L}=\mu_n(H)$.
Recall from \cite{MN04} that
\begin{equation}\label{eq:qpLambda}
\begin{aligned}
\Delta(\Lambda)p_Y(a \o 1) &= \Delta(\Lambda)p_Y(1 \o S(a)), \\
 (1 \o a)q_X\Delta(\Lambda)  & = (S(a)\o 1)q_X\Delta(\Lambda)
\end{aligned}
\end{equation}
for all $a \in H$. Then we have
\begin{eqnarray*}
T_{L,L}&=& m\left((\id\otimes \Delta^{(n-1)})(q_L)\Delta^{(n)}(\Lambda)\phi_n(\Delta^{(n-1)}\o \id)(p_L)\right)\\
&=& m\left((\id\otimes \Delta^{(n-1)})(q_L)\phi_n(\Delta^{(n-1)}\o \id)(\Delta(\Lambda)p_L)\right)
\quad
\mbox{by \eqref{eq:phi_nLambda}}\\
& =& m\left((q_L^{(1)}\otimes 1 )\phi_n(\Delta^{(n-1)}\o \id)(\Delta(\Lambda)p_L (q_L^{(2)}\o 1))\right)\quad
\mbox{by Lemma \ref{product_trick}}\\
& =& m\left((q_L^{(1)}\otimes 1 )\phi_n(\Delta^{(n-1)}\o \id)(\Delta(\Lambda)p_L (1 \o S(q_L^{(2)})))\right)
\quad
\mbox{by \eqref{eq:qpLambda}}\\
& =& q_L^{(1)}m\left(\Delta^{(n)}(\Lambda)\phi_n(\Delta^{(n-1)}\o \id)(p_L) \right)S(q_L^{(2)})\\
& =& q_L^{(1)}m\left((\id \o \Delta^{(n-1)})((\b q_R^{(1)} \o q_R^{(2)})
\Delta(\Lambda))\phi_n(\Delta^{(n-1)}\o \id)(p_L) \right)S(q_L^{(2)})\quad
\mbox{by \eqref{eq:q_Lambda_p}}\\
&= & q_L^{(1)}\b \mu_n(H) S(q_L^{(2)})\\
&= & q_L^{(1)}\b  S(q_L^{(2)})\mu_n(H)\quad\mbox{by Theorem \ref{t:1}}\\
& = & \mu_n(H)\,.
\end{eqnarray*}
Thus, we have $T_{X,L}=\mu_n(H)$ for $X \in \{L, R\}$. Now, for $Y=L$ or $R$, we have
\begin{eqnarray*}
T_{X,Y}&=& m\left((\id\otimes \Delta^{(n-1)})(q_X)\Delta^{(n)}(\Lambda)\phi_n(\Delta^{(n-1)}\o \id)(p_Y)\right)\\
&=& m\left((\id\otimes \Delta^{(n-1)})(q_X)\Delta^{(n)}(\Lambda)\phi_n (1 \o p_Y^{(2)})
(\Delta^{(n-1)}(p_Y^{(1)}) \o 1)\right)\\
&=& m\left((\id\otimes \Delta^{(n-1)})((1 \o p_Y^{(1)})q_X)\Delta^{(n)}(\Lambda)\phi_n (1 \o p_Y^{(2)})
\right)\quad
\mbox{by Lemma \ref{product_trick}}\\
&=& m\left((\id\otimes \Delta^{(n-1)})((1 \o p_Y^{(1)})q_X\Delta(\Lambda))\phi_n (1 \o p_Y^{(2)})
\right)\\
&=& m\left((\id\otimes \Delta^{(n-1)})(( S(p_Y^{(1)})\o 1)q_X\Delta(\Lambda))\phi_n (1 \o p_Y^{(2)})
\right)\quad
\mbox{by \eqref{eq:qpLambda}}\\
&=&  S(p_Y^{(1)}) m\left((\id\otimes \Delta^{(n-1)})(q_X\Delta(\Lambda))\phi_n
\right) p_Y^{(2)}\\
&=&  S(p_Y^{(1)}) m\left((\id\otimes \Delta^{(n-1)})(q_X)\Delta^{(n)}(\Lambda)\phi_n
\right) p_Y^{(2)}\\
&=&  S(p_Y^{(1)}) m\left((\id\otimes \Delta^{(n-1)})(q_X)\phi_n(\Delta^{(n-1)}\o \id)
(\Delta(\Lambda)(p_L^{(1)} \o p_L^{(2)}\a))\right) p_Y^{(2)}
\quad \mbox{by \eqref{eq:q_Lambda_p}} \\
&=&  S(p_Y^{(1)}) \mu_n(H)\a  p_Y^{(2)}\\
&=& \mu_n(H)
\end{eqnarray*}
This proves the statement.
\end{proof}

\section{Hopf algebras with central group-like elements}\label{s:HCG}
Let $H$ be a finite-dimensional Hopf algebra over $\BC$ and $G(H)$
the group of all group-like elements of $H$. Let $G$ be a subgroup
of $G(H)$ which lies in the center of $H$, $\w$ a normalized
3-cocycle of $G$ with coefficients in $\BC^\times$ and $j:G
\rightarrow \hat{G}$ a group isomorphism. Then $j$ can be extended
to a Hopf algebra isomorphism from $\BC[G]$ to
$\BC[\hat{G}]=\BC[G]^*$. Let
$$
e_x = \frac{1}{|G|}\sum_{y \in G} j(y)(x)^{-1} y\,.
$$
By the orthogonality of characters of finite groups,
$$
j(e_x)(z) = \frac{1}{|G|}\sum_{y \in G} j(y)(x)^{-1} j(y)(z) =
\frac{1}{|G|}\sum_{y \in G} j(y)(x^{-1}z) =\delta_{x,z}\,.
$$
Therefore, $\{j(e_x)| x \in G\}$ is the dual basis of $G$ for
$\BC[G]^*$. In particular, $\{e_x\}_{x \in G}$ is the complete set
of orthogonal primitive idempotents of $\BC[G]$ and we have the
equalities
\begin{equation}\label{eq:properties}
S(e_x)=e_{x\inv}, \quad \e(e_x)=\delta_{1,x}, \quad \Delta(e_x) = \sum_{y \in G} e_{xy\inv} \o e_y
\end{equation}
for all $x \in G$. One can construct a new quasi-Hopf algebra using these data of $H$. The following lemma is
a joint observation with Geoffrey Mason.

\begin{lem}
  Let $H$ be a Hopf algebra over $\BC$ and $G$ a subgroup of $G(H)$. If $G$ lies in the
  center of $H$, then for any normalized 3-cocycle $\w$ of $G$ with coefficient in $\BC^\times$
  and for any group isomorphism $j: G \rightarrow \hat{G}$, the tuple $H_{(G, \w, j)}=(H, \Delta, \e, \phi, \a, \b, S)$ is a
  quasi-Hopf algebra where $\phi$, $\a$ and $\b$ are defined by
 \begin{equation}\label{eq:phi}
\phi=\sum_{x, y, z \in G}\w(x, y, z)\inv e_x \o e_y \o e_z, \quad
\a=1,\quad \b = \sum_{x \in G}\w(x, x\inv, x) e_x\,,
\end{equation}
   and $\Delta$, $\e$ and $S$
  are the comultiplication, counit and antipode of $H$. Moreover, if $\w'$ is
  cohomologous to $\w$, then $H_{(G, \w, j)}$ and $H_{(G, \w', j)}$ are gauge equivalent quasi-Hopf
  algebras.
\end{lem}
\begin{proof}
Note that $\phi\inv =\sum_{x, y, z \in G}\w(x, y, z) e_x \o e_y \o
e_z$ and $\b\inv= \sum_{x \in G}\w(x, x\inv, x)\inv e_x$. Using the
fact that $\BC[G]$ lies in the center of $H$ and
\eqref{eq:properties}, it is straightforward to verify that $H_{(G,
\w, j)}$ is a quasi-Hopf algebra. Let $\w'=\w \delta b$ for some
normalized 2-cochain  $b$ of $G$ with coefficients in $\BC^\times$
where
$$
\delta b (x,y,z)= \frac{b(y,z)b(x,yz)}{b(xy,z)b(x,y)}
$$
for $x, y, z \in G$. We define
$$
F=\sum_{x,y \in G} b(x,y)\inv e_x \o e_y\,.
$$
Since $b$ is a normalized cochain, $F$ is a gauge transformation on
$H_{(G, \w, j)}$. Moreover,
$$
\phi^F=\sum_{x,y,z \in G} \w(x,y,z)\inv \delta b\inv(x,y,z)\, e_x \o e_y \o e_z =
\sum_{x,y,z \in G} \w'(x,y,z)\inv \, e_x \o e_y \o e_z\,.
$$
Since $F$ lies in the center of $H \o H$, $\Delta^F=\Delta$ and so
$H_{(G, \w',j)} = H_{(G, \w,j)}^F$  as quasi-bialgebras.
\end{proof}
\begin{remark}
  $H$ and $H_{(G, \w, j)}$ are identical as bialgebras and so their module categories are
  linearly
  equivalent. If $H$ is also semisimple, then $H$ and $H_{(G, \w, j)}$ have the same fusion rules for their
  irreducible representations. However, in general, $\C{H}$ and $\C{H_{(G, \w, j)}}$ are not equivalent
  tensor categories if $\w$ is not a coboundary. We will see
  examples for this below.
\end{remark}

Let us further assume that $H$ is a semisimple Hopf algebra over
$\BC$ with a central group-like element $u$ of order 2 and  $G$ is
the subgroup generated by $u$. Note that there is exactly one group
isomorphism $j$ from $G$ to
$\hat{G}$ and $H^3(G, \BC^\times)$ is an abelian group of order 2.\\

Consider the function $\w: G \times G \times G \rightarrow \BC^\times$ defined by
$$
\w(x, y, z)=\left\{\begin{array}{rl}
  -1 & \mbox{ if } x=y=z=u,\\
  1 & \mbox{ otherwise}\,.
\end{array}\right.
$$
One can easily verify that $\w$ is a non-trivial 3-cocycle of $G$
and so $\w$ represents the unique non-trivial cohomology class of
$H^3(G, \BC^\times)$. We will simply write $H_u$ for the quasi-Hopf
algebra $H_{(G,\w, j)}$.\\

Let $\chi$ be the non-trivial character of $G$.
Then $j(g)=\chi$ and so
\begin{equation}\label{eq:beta}
e_1 = \frac{1}{2}(1+u), \quad e_g=\frac{1}{2}(1-u) \quad \mbox{and}
\quad \b=e_1-e_u = u\,.
\end{equation}

We will proceed to obtain a formula of $\mu_n(H_u)$. Let us denote
the $n$-th Sweedler power of an element $x$ of $H$ by $x^{[n]}$. It
is easy to see that
$$
x^{[n]}=\left\{\begin{array}{cl}
  x & \mbox{ if } $n$ \mbox{ is odd,}\\
  \e(x)1 & \mbox{ if } $n$ \mbox{ is even}
\end{array}\right.
$$
for any $x \in \BC[G]$. In particular,
$$
e_z^{[n]}=\left\{\begin{array}{cl}
  e_z & \mbox{ if } $n$ \mbox{ is odd,}\\
  \delta_{z,1}1 & \mbox{ if } $n$ \mbox{ is even}
\end{array}\right.
$$
for any $z \in G$. Thus, for any positive integer $n$, we have
$$
m\left(\phi^{(1)}\o \Delta^{(n)}(\phi^{(2)})\o \phi^{(3)}\right) =
\sum_{x,y,z \in G} \w(x, y, z)\inv e_x e_y^{[n]} e_z = u^{n}\,.
$$
where $m$ is the multiplication of $H$. Since $\BC[G]$ is
commutative, by \eqref{eq:phi_n},  we have
\begin{equation}\label{eq:mphi_n}
m(\phi_r)=u^{(r-1)(r-2)/2}
\end{equation}
for all $r\geq 1$. Indeed the formula is clearly correct for $r=1,
2$, and for $r\geq 2$ we have inductively
\begin{align*}
  m(\phi_{r+1})&=m\bigl((1\o\phi_r)(\phi\so
  1\o\Delta^{(r-1)}(\phi\so2)\o\phi\so3)\bigr)\\
  &=m(\phi_r)m\bigl(\phi\so
  1\o\Delta^{(r-1)}(\phi\so2)\o\phi\so3\bigr)\\
  &=u^{(r-1)(r-2)/2}u^{r-1}=u^{r(r-1)/2}\,.
\end{align*}

\begin{prop}\label{p1}
  Let $H$ be a finite-dimensional semisimple Hopf algebra over $\BC$ with a central group-like element $u$ of order 2.
  Suppose that $V$ is a finite-dimensional simple $H$-module with character $\chi$,
  and that ${\overline V}$ is
  the $H_u$-module associated with $V$.
   Then the $n$-th Frobenius-Schur indicator of
  ${\overline V}$ is given by
  $$
  \nu_n({\overline V}) = \nu_n(V)\chi(u^{(n-3)n/2})/\chi(1)
  $$
  for any positive integer $n$, where $\nu_n(V)$ is the $n$-th
  Frobenius-Schur indicator of $V$ considered as an $H$-module.
  \end{prop}
  \begin{proof}
   Since $u$ is in the center of $H$, $\phi_n$ is also in the center of $H^{\o n}$. Thus, by Theorem
  \ref{t:1} and \eqref{eq:mphi_n}, we have
  $$
  \begin{aligned}
  \mu_n(H_u)&= m(\Delta^{(n)}(\Lambda)\phi_n)(\b \a)\inv = m(\Delta^{(n)}(\Lambda))m(\phi_n) u^{-1} \\
   &= \Lambda^{[n]} u^{(n-2)(n-1)/2} u\inv = \Lambda^{[n]}u^{(n-3)n/2}\,,
  \end{aligned}
  $$
  where $\Lambda$ is the normalized integral of $H$. The second
  and the third equalities follow from \eqref{eq:beta} and \eqref{eq:mphi_n}.
  By Theorem \ref{t:1},
  $$
  \begin{aligned}
  \nu_n({\overline  V}) &= \chi(\Lambda^{[n]}u^{(n-3)n/2)})\\
                      & = \chi(\Lambda^{[n]})\chi(u^{(n-3)n/2)}/\chi(1) \\
                      &= \nu_n(V)\chi(u^{(n-3)n/2})/\chi(1)\,.
  \end{aligned}
  $$
  The second equality follows from the fact that $u$ acts on $V$ as
  a scalar.
\end{proof}

\begin{example}
As the simplest example, consider $H=\BC[G]$, where $G$ is the group
of order 2 generated by $u$. Let $V$ be the non-trivial
1-dimensional $H$-module. Since $u$ acts on $V$ as the scalar $-1$,
we have $\nu_n(V)=\frac{1+(-1)^n}{2}$  for all positive integer $n$.
By Proposition \ref{p1}, we have
$$
\nu_n(\overline
V)=\frac{1+(-1)^n}{2}(-1)^{n(n-3)/2}=\cos\left(\frac{n\pi}{2}\right)\,.
$$
\end{example}

\section{Examples and Applications}\label{ea}
In \cite{TaYa98}, Tambara and Yamagami has classified that there are
four inequivalent fusion categories over $\BC$ with five simple
objects $\{a, b, c, d,  m\}$ and fusion rules:
  \begin{enumerate}
    \item $am=ma=bm=mb=cm=mc=dm=md=m$,
    \item $mm=a+b+c+d$, and
    \item $\{a, b, c, d\}$ forms a abelian group isomorphic to $\BZ_2 \times \BZ_2$\,.
  \end{enumerate}
Three of these categories are tensor equivalent to the module
categories of the following Hopf algebras: $\BC[Q_8]$, $\BC[D_8]$
and the 8-dimensional Kac algebra $K$, where $Q_8$ and $D_8$ are,
respectively, the quaternion group and the dihedral group of order
8. It was a question raised by Susan Montgomery whether one can
explicitly construct a quasi-Hopf algebra whose module category is
tensor equivalent to the fourth fusion category with above fusion
rules. In this section, we will answer this question by showing that
these four fusion categories are tensor equivalent to
$$
\C{\BC[D_8]}, \quad  \C{\BC[Q_8]}, \quad \C{K}, \quad \mbox{and}
\quad \C{K_u},
$$
where $u$ is the unique order 2 central  group-like element of $K$.
In particular, we show that $\BC[D_8]$ is gauge equivalent to the
quasi-Hopf algebra $\BC[Q_8]_u$ and $\BC[D_8]_u$ is gauge
equivalent to $\BC[Q_8]$\\

Let $H=\BC[Q_8]$, $\BC[D_8]$ or $K$. Then $H$ has a unique central
group-like element $u$ of order 2 and it has only one 2-dimensional
irreducible representation $V$.
We proceed to compute the Frobenius-Schur indicators of $V$ in $\C{H}$ and
 those of ${\overline V}$ in $\C{H_u}$.\\

The 8-dimensional Kac algebra $K$ is a semisimple Hopf algebra over $\BC$
generated by $x, y, z$ as $\BC$-algebra with the relations:
$$
x^2=1, \quad y^2=1, \quad z^2=\frac{1}{2}(1 + x+y -xy), \quad xy=yx, \quad xz=zy, \quad yz=zx\,.
$$
The coalgebra structure is given by
$$
\Delta(x)=x \otimes x, \quad \Delta(y)=y \otimes y, \quad \e(x)=\e(y)=1\,,
$$
$$
\Delta(z)=\frac{1}{2}(1 \o 1 + x \o 1 + 1 \o y - y \o x)(z \o z), \quad \e(z)=1\,.
$$
The antipode of $K$ is determined by
$$
S(x)=x, \quad S(y)=y, \quad S(z)=z\,.
$$
Note that $u=xy$ is the unique central group-like element of order 2
of $K$ (cf. \cite{Mas95} and \cite{TaYa98} for more details on this
Kac algebra). Let $V$ be the degree 2 irreducible representation of
$K$ and $\chi$ the character of $V$. We have $\chi(u)=-2$. The
higher Frobenius-Schur indicators of $V$ considered as a $K$-module
are given by
$$
\chi(\Lambda^{[n]}) = 1, 0, 0, 0, 1, 0, 2
$$
where $n=2,\dots, 8$ and $\Lambda$ is the normalized integral of $K$
given by
$$
\Lambda= \frac{1}{8} ((1+x+y+xy)+ (1+x+y+xy)z)\,.
$$
By Proposition \ref{p1}, the $n$-th Frobenius-Schur indicators
$(n=2,\dots,8)$ of $V$ considered as a $K_u$-module are
$$
-1, 0, 0, 0, -1, 0, 2\,.
$$
\\
Since $Q_8$ and $D_8$ have centers of order 2, each of $\BC[Q_8]$
and $\BC[D_8]$ has a unique central group-like element $u$ of order
2. Let $H=\BC[Q_8]$ or $\BC[D_8]$,  $V$ the degree 2 irreducible
representation of $H$ and $\chi$ the character of $V$. Then $\chi(u)=-2$.\\

For $D_8$, the Frobenius-Schur indicators $\nu_n(V)$ of
$V$($n=2,\dots, 8$) are well-known to be
$$
1, 0, 2, 0, 1, 0, 2\,.
$$
The Frobenius-Schur indicators $\nu_n(V)$ of $V$ ($n=2,\dots, 8$)
for $Q_8$ are
$$
-1, 0, 2, 0, -1, 0, 2\,.
$$

By Proposition \ref{p1}, we can complete the following table of the
Frobenius-Schur indicators for the 2-dimensional irreducible
representation $V$ of these quasi-Hopf algebras.

$$
\begin{array}{|c|c|c|c|c|c|c|c|}
\hline
 & \nu_2(V)& \nu_3(V) & \nu_4(V)& \nu_5(V) & \nu_6(V)& \nu_7(V)& \nu_8(V)\\
 \hline
 K & 1 &  0 & 0 & 0 & 1 & 0 & 2\\
  \hline
 K_u & -1 &  0 & 0 & 0 & -1 & 0 & 2\\
  \hline
\BC[D_8] & 1& 0& 2& 0& 1& 0& 2\\
\hline
\BC[D_8]_u & -1& 0& 2& 0& -1& 0& 2\\
\hline
\BC[Q_8] & -1& 0& 2& 0& -1& 0& 2\\
\hline
\BC[Q_8]_u & 1& 0& 2& 0& 1& 0& 2\\
\hline
\end{array}
$$
\newline
As a conclusion from this table, we have the following theorem.
\begin{thm}
 The fusion categories
$$
\C{K},\quad \C{K_u},\quad \C{\BC[D_8]},\quad \C{\BC[Q_8]}
$$
form a complete set of inequivalent fusion categories with the
fusion rules given at the beginning of this section. Moreover, the
quasi-Hopf algebra $\BC[D_8]_g$ is gauge equivalent to  $\BC[Q_8]$
and $\BC[Q_8]_g$ is gauge equivalent to  $\BC[D_8]$
\end{thm}
\begin{proof}
It follows from the above table and Proposition \ref{p:invarianceI} that the fusion categories
$$
\C{K},\quad \C{K_u},\quad \C{\BC[D_8]},\quad \C{\BC[Q_8]}
$$
are inequivalent tensor categories. Since there are totally four
inequivalent fusion categories with the same fusion rules given at
the beginning of this section (cf. \cite{TaYa98}), the categories
$\C{K}$, $\C{K_u}$, $\C{\BC[D_8]}$, $\C{\BC[D_8]_u}$ account for all
these fusion categories. Since the Frobenius-Schur indicators of the
degree 2 irreducible representation of $\BC[Q_8]_u$ is the same as
that of $\BC[D_8]$ and the fusion rules of $\C{\BC[Q_8]_u}$ and
$\C{\BC[D_8]}$ are the same, the categories $\C{\BC[D_8]}$ and
$\C{\BC[Q_8]_u}$ are tensor equivalent. By Theorem
\ref{t:gauge_equivalence}, $\BC[D_8]$ and $\BC[Q_8]_u$ are gauge
equivalent quasi-Hopf algebras. Similarly, one can show that the
quasi-Hopf algebras $\BC[Q_8]$ and $\BC[D_8]_u$ are also gauge
equivalent.
\end{proof}

\section{Higher Indicators for quasi-Hopf algebras associated to group cocycles}\label{s:DG}
 In this section we will derive explicit Frobenius-Schur
(FS)-indicator formulae for two kinds of quasi-Hopf algebras
associated to a 3-cocycle on a finite group $G$. In particular,
Bantay's formula of the second FS-indicator for a twisted quantum
double of a finite group derived in \cite{MN04} is a special case of
these formulae (see also
 \cite{Bantay00} and \cite{Bantay97}).\\

 Let $G$ be a finite group and
$\w$ a normalized 3-cocycle on $G$ with coefficients in
$\BC^\times$. One can construct a quasi-Hopf algebra
$H(G,\w)=(\BC[G]^*, \Delta, \e, \phi, \a, \b, S)$ where
multiplication, identity, comultiplication, counit and antipode
are the same as the structure maps of the Hopf algebra $\BC[G]^*$,
and $\phi$, $\a$, and $\b$ are given by
$$
\phi=\sum_{a,b,c \in G}\w(a,b,c)e(a) \o e(b) \o e(c)\, \quad \a=1,
\quad\mbox{and} \quad \b=\sum_{a\in G}\w(a, a\inv, a)\inv e(a)\,,
$$
where $\{e(x)\mid x \in G\}$ is the dual basis of $G$ for
$\BC[G]^*$. The quantum double of $H(G,\w)$ is called the twisted
quantum double $D^{\w\inv}(G)$.

\subsection*{The quasi-Hopf algebras $H(G,\w)$}
The basis $\{e(g)\mid g \in G\}$ for $H(G,\w)$ is a complete set of
primitive idempotents of $H(G,\w)$ and its dual basis $\{\chi_x \mid
x \in G\}$ is the complete set of irreducible characters of
$H(G,\w)$.
\begin{prop}
  Let $G$ be a finite group and $\w$ a normalized 3-cocycle on $G$ with coefficients in $\BC^\times$.
  For any simple
  $H(G,\w)$-module $V_x$ with character $\chi_x$ and positive integer $n$,
  \begin{equation}\label{eq:FSI_formula}
  \nu_n(V_x)=\delta_{x^n,1} \prod_{r=1}^{n-1} \w(x, x^r, x)\,.
  \end{equation}
  In particular, if $n$ is not a multiple of the order $o(x)$ of $x$,
  then $\nu_n(V_x)$ vanishes. $\nu_{o(x)}(V_x)$ is a root of
  unity whose order is the same as that of the cohomology class
  $\res_{\langle x \rangle}[\w]\in H^3(\langle x\rangle,\BC^\times)$, and
  $\nu_{s\cdot o(x)}(V_x)=\nu_{o(x)}(V_x)^s$ for $s\in\mathbb N$.
\end{prop}
\begin{proof}
Clearly, $e(1)$ is the normalized integral of $H(G,\w)$. For any $x
\in G$ and integer $n \ge 1$,
$$
\e(x)^{[n]}=m\left(\Delta^{(n)}(e(x))\right)=\sum_{y^n=x} e(y)\,,
$$
where $m$ denotes the multiplication of $H(G,\w)$,
$\Delta^{(1)}=\id$ and $\Delta^{(n+1)}=(\id^{\o n-1}\o
\Delta)\Delta^{(n)}$ for $n \ge 1$. By the commutativity of
$H(G,\w)$ and Theorem \ref{t:1}, we have
$$
\nu_n(V_x) = \chi_x\left( e(1)^{[n]} \b\inv  m(\phi_n)\right)\,,
$$
where
$$
\phi_1=1_H, \quad \phi_2=1 \o 1, \quad \phi_{r+1}=(1 \o \phi_n)(\id
\o \Delta^{(r-1)}\o \id)(\phi)\,.
$$
for any integer $r \ge 2$. Note that
$$
m\left((\id \o \Delta^{(r)}\o \id)(\phi)\right)=\sum_{a, b, c\in G}
\w(a,b,c)e(a)e(b)^{[r]}e(c) =\sum_{y\in G} \w(y,y^r,y)e(y)\,.
$$
Hence, by induction,
$$
m(\phi_n)=\sum_{y \in G} \prod_{r=1}^{n-2}\w(y, y^r, y) e(y)\,.
$$
Thus, we have
$$
\nu_n(V_x)=\chi_x\left(\sum_{y^n=1} e(y) \w(y,y\inv, y)
\prod_{r=1}^{n-2}\w(y, y^r, y)\right)=
\delta_{x^n,1}\prod_{r=1}^{n-1}\w(x, x^r, x)\,.
$$
Following the description in \cite{MS89}, the class of the 3-cocycle
$\w_x$ of the cyclic subgroup $\langle x \rangle$ of order $N$
defined by
$$
\w_x(x^\ell, x^m, x^n)=\exp\left(\frac{2\pi i}{N^2}\hh\ell(\hh m+\hh
n-\hh{m+n})\right)
$$
generates the group $H^3(\langle x \rangle,\BC^\times)$, where $\hh
n$ denotes the remainder upon the division of $n$ by  $N$. Then we
have
$$
\w_x(x, x^r,x)=\exp\left(\frac{2\pi  i}{N^2}(\hh
r+1-\hh{r+1})\right)
 =\begin{cases}
 \exp\left(\frac{2\pi i}{N}\right)=\zeta &\text{ if } r\equiv -1\mod N,\\
 1&\text{ otherwise.}
 \end{cases}
$$
Obviously, there is a 3-coboundary $f$ of $G$ such that
$$
\w f=\w_x^t \quad \mbox{on $\langle x \rangle$ for some integer
$t$.}
$$
Since $\nu_n(V_x)$ is a gauge invariant,
$$
\nu_{Ns}(V_x)= \prod_{r=1}^{Ns-1}\w_x^t(x, x^r, x) = (\zeta^t)^s.
$$
Of course, the order of $\zeta^t$ is the same as the order of
$\res_{\langle x \rangle}[\w]$.
\end{proof}
\subsection*{The twisted quantum double of finite group $D^\w(G)$}
The {\em twisted quantum double} $D^\w(G)$ of $G$ with respect to
$\w$ is the semisimple quasi-Hopf algebra with underlying vector
space $\BC[G]^* \o \BC[G]$  in which multiplication,
comultiplication $\Delta$, associator $\phi$, counit $\e$, antipode
$S$, $\a$ and $\b$ are given by
\begin{equation}\label{multiplication}
(e(g) \o x)(e(h) \o y) =\theta_g(x,y) \delta_{g^x,h}
e(g)\o x y\,,
\end{equation}
\begin{equation}\label{comultiplication}
\Delta(e(g)\o x)  = \sum_{hk=g} \g_x(h,k) e(h)\o x
\o e(k) \o x\,,
\end{equation}
\begin{equation}\label{eq0.01}
\phi = \sum_{g,h,k \in G} \w(g,h,k)^{-1} e(g) \o 1 \o
e(h) \o 1 \o e(k) \o 1\,,
\end{equation}
\begin{equation}\label{eq:antipode}
\begin{aligned}
& \e(e(g)\o x) = \delta_{g,1}, \quad \a=1, \quad \b=\sum_{g \in
G} \w(g,g^{-1},g)e(g)\o 1\,,\\
& S(e(g)\o x) =
\theta_{g^{-1}}(x,x^{-1})^{-1}\g_x(g,g^{-1})^{-1}e(x^{-1}g^{-1}x)\o
x^{-1}\,,
\end{aligned}
\end{equation}
where $\delta_{g,1}$ is the Kronecker delta,  $g^x=x\inv g x$, and
\begin{eqnarray*}
\theta_g(x,y) &=&\frac{\w(g,x,y)\w(x,y,(x y)^{-1}g x y)}{\w(x,x^{-1}g x,y)}, \\
\g_g(x,y) & = & \frac{\w(x,y,g)\w(g, g^{-1}x g, g^{-1}yg)}{\w(x,g,
g^{-1}y g)}
\end{eqnarray*}
for any $x, y, g \in G$ (cf. \cite{DPR91}).\\

By induction, one can show that
\begin{equation}\label{eq:DG_phi_n}
\phi_n= \sum_{a_1, \cdots, a_{n} \in G} \left(\prod_{i=1}^{n-2}
\w(a_i,a_{i+1}\cdots a_{n-1},a_{n})^{-1}\right)e(a_1) \o 1
\o \cdots \o e(a_{n}) \o 1
\end{equation}
for any integer $n \ge 2$.
The normalized integral of $D^\w(G)$ is given by
$$
\Lambda = \frac{1}{|G|}\sum_{x \in G} e(1) \o x\,.
$$
Therefore,
$$
\Delta(\Lambda)= \frac{1}{|G|}\sum_{x,a \in G} \g_x(a\inv, a)e(a\inv)
\o x \o e(a) \o x
$$
and, by induction, we have
\begin{equation}\label{eq:DG_Delta_n}
\begin{aligned}
  \Delta^{(n)}(\Lambda) =
   \frac{1}{|G|}\sum_{{x, a_1,
\dots, a_n \in G} \atop {1=a_1\cdots a_n}} \left(\prod_{i=1}^{n-1}
\g_x(a_i, a_{i+1}\cdots a_n)\right) e(a_1) \o x  \cdots \o e(a_n) \o
x
\end{aligned}
\end{equation}
for any integer $n \ge 1$.
Thus, we have
\begin{prop}
Let $G$ be a finite group and $\w$ a normalized 3-cocycle of $G$
with coefficients in $\BC^\times$. Then, for any integer $n \ge 2$,
the $n$-th Frobenius-Schur indicator of any representation $V$ of
$D^\w(G)$ with character $\chi$ is given by
\begin{equation}\label{eq:DG_FSI}
\begin{aligned}
\nu_n(V) = \frac{1}{|G|}\sum_{x, a \in G \atop{x^{-n}=(ax\inv)^n}} &
\left( \prod_{i=1}^{n-2} \frac{\g_x(a^{x^i}, a^{x^{i+1}}\cdots
a^{x^{n-1}})\theta_a(x^i,x)} {\w(a^{x^{i}}, a^{x^{i+1}}\cdots
a^{x^{n-1}},a^{x^n})}\right)\\
& \qquad \cdot \frac{\g_x(a ,
a\inv)\theta_a(x^{n-1},x)}{\w(a, a\inv, a)}
\cdot \chi\left(e(a) \o x^n\right)\,.
\end{aligned}
\end{equation}
\end{prop}
\begin{proof}
  Since both $\a$ and $\b$ of $D^\w(G)$ are invertible, by Theorem \ref{t:1}, we have
\begin{equation}\label{eq:simplified}
  \mu_n(D^\w(G))= \b\inv m(\Delta^{(n)}(\Lambda) \phi_n)\,.
\end{equation}
By \eqref{eq:DG_phi_n} and \eqref{eq:DG_Delta_n}, we obtain that
$$
\begin{aligned}
\Delta^{(n)}(\Lambda) \phi_n=
  \frac{1}{|G|}\sum_{x, a_1,
\cdots a_n \in G \atop{1= a_1 \cdots a_{n}}}& \left(
\prod_{i=1}^{n-2} \frac{\g_x(a_{i+1}, a_{i+2}\cdots a_n)}
{\w(a_i^x ,a_{i+1}^x\cdots a_{n-1}^x,a_n^x)}\right)\cdot \g_x(a_1 , a_1\inv) \\
& \qquad \cdot e(a_1) \o x \o \cdots \o e(a_{n}) \o x\,.
\end{aligned}
$$
Thus we have
$$
\begin{aligned}[b]
\b\inv m\left(\Delta^{(n)}(\Lambda) \phi_n\right)& =
\begin{aligned}[t]
 \frac{1}{|G|}\!\sum\limits_{x, a \in G \atop{aa^x\cdot a^{x^{n-1}}=1}}
 & \left(\prod_{i=1}^{n-2} \frac{\g_x(a^{x^i}, a^{x^{i+1}}\cdots
a^{x^{n-1}})\theta_a(x^i,x)} {\w(a^{x^{i}}, a^{x^{i+1}}\cdots
a^{x^{n-1}},a^{x^n})}\right)\\
& \cdot \g_x(a , a\inv)\theta_a(x^{n-1},x)
\cdot \b\inv \cdot (e(a) \o x^n)
\end{aligned}\\ \\
& = \begin{aligned}[t]
 \frac{1}{|G|}\sum\limits_{x, a \in G \atop{x^{-n}=(ax\inv)^n}}
  \biggl(\prod_{i=1}^{n-2}  & \frac{\g_x(a^{x^i}, a^{x^{i+1}}\cdots
a^{x^{n-1}})\theta_a(x^i,x)} {\w(a^{x^{i}}, a^{x^{i+1}}\cdots
a^{x^{n-1}},a^{x^n})}\biggr)\\
& \cdot \frac{\g_x(a , a\inv)\theta_a(x^{n-1},x)}{\w(a,a\inv, a)}
 \cdot (e(a) \o x^n)
 \end{aligned}
\end{aligned}
$$
The statement follows easily  from \eqref{eq:simplified} and Theorem \ref{t:1}.
\end{proof}
\begin{remark}
  The formula for the higher indicators may look different if one uses another form of $\mu_n(H)$.
  For example, $\mu_n(D^\w(G)) = \b\inv m(\phi_n (\Delta^{(n-1)}\o \id)\Delta(\Lambda))$ since
  $\phi_n (\Delta^{(n-1)}\o \id)\Delta(\Lambda)= \Delta^{(n)}(\Lambda)\phi_n$. Using this form of the cental
  invariant $\mu_n(D^\w(G))$, one will obtain
  \begin{equation}
\begin{aligned}
\nu_n(V) = \frac{1}{|G|}\sum_{(ax^{-1})^{n}=x^{-n}}
 & \left( \prod_{i=1}^{n-2}
\frac{\g_x(a^{x^{i-1}},a^{x^i}\cdots a^{x^{n-2}})\theta_a(x^i,x)}
{\w(a^{x^{i-1}},a^{x^i}\cdots a^{x^{n-2}},a^{x^{n-1}})}
\right)\\
&  \cdot \frac{\g_x(a^{-x^{n-1}} , a^{x^{n-1}})\theta_a(x^{n-1},x)}{\w(a, a\inv, a)}
\cdot \chi\left(e(a) \o x^n\right)\,.
\end{aligned}
\end{equation}
\end{remark}

%GATHER{MYBIBL.BIB}
\bibliographystyle{amsalpha}
%\bibliography{mybibl}
\providecommand{\bysame}{\leavevmode\hbox
to3em{\hrulefill}\thinspace}
\providecommand{\MR}{\relax\ifhmode\unskip\space\fi MR }
\providecommand{\MRhref}[2]{%
  \href{http://www.ams.org/mathscinet-getitem?mr=#1}{#2}
} \providecommand{\href}[2]{#2}

\end{document}